\documentclass[12pt]{article}
\usepackage{latexsym}
\usepackage{amsmath,epsfig}
\usepackage{amssymb}
\usepackage{amsfonts}

\numberwithin{equation}{section}

\newtheorem{theorem}{Theorem}[section]
\newtheorem{lemma}{Lemma}[section]

\newcommand{\eproof}{{\mbox{\ }~\hfill
\mbox{\large $\Box$} \par \vskip 10pt}}

\title{Increasing stability for determining the potential in the Schr\"odinger equation with attenuation from the Dirichlet-to-Neumann map}
\author{ Victor Isakov\thanks{Department of Mathematics and Statistics, Wichita State University, KS 67260-0033, USA. Email:victor.isakov@wichita.edu}\qquad Jenn-Nan Wang\thanks{Institute of Applied Mathematics, NCTS (Taipei), National
Taiwan University, Taipei 106, Taiwan. Email:jnwang@math.ntu.edu.tw}}

\date{}
\begin{document}
\renewcommand{\theequation}{\thesection.\arabic{equation}}
 
\maketitle
\begin{abstract}
We derive some bounds which can be viewed as an evidence of increasing stability in the problem of recovering  the potential coefficient in the Schr\"odinger equation from the Dirichlet-to-Neumann map in the presence of attenuation, when energy level/frequency is growing. These bounds 
hold  under certain a-priori regularity constraints on the unknown coefficient. Proofs use complex and bounded complex geometrical optics solutions. 
\end{abstract}

\section{Introduction}\label{sec1}
\setcounter{equation}{0}

We consider the problem of recovery of the potential in the Schr\"odinger equation with attenuation from the Dirichlet-to-Neumann map.  The use of complex exponential solutions plays essential role in the development of this problem. The idea is dated back to Calderon and Faddeev who introduced  complex exponential solutions to demonstrate uniqueness in the closely related linearized inverse conductivity problem and in the inverse potential scattering problem for the Schr\"odinger equation.  A breakthrough was made by Sylvester and Uhlmann in \cite{SU} where they constructed almost complex exponential solutions, and proved global uniqueness of $c$ (potential  in the Schr\"odinger equation) in the three-dimensional case.  A logarithmic stability estimate for $c$ from the Dirichlet-to-Neumann map was obtained by Alessandrini \cite{A} and the optimality of log-type stability (at zero energy) was demonstrated by Mandache \cite{M}. The logarithmic stability is quite discouraging for applications, since small errors in the data of the inverse problem result in  large errors in numerical reconstruction of physical properties of the medium. In particular, it severely restricts resolution in the electrical impedance tomography.

For the problem of recovering the potential in the Schr\"odinger equation without attenuation from the Dirichlet-to-Neumann map, the first author in \cite{I3} derived some stability estimates in different ranges of frequency, which demonstrate the increasing stability phenomena as the frequency/energy $k$ is growing. In \cite{I3}, both complex- and real-valued geometrical optics solutions were used in the proof. The proof was simplified in \cite{INUW} where only complex-valued geometrical optics solutions were used. Similar results were obtained by Isaev and Novikov \cite{INo} by  less explicit and more complicated methods of scattering theory.

Continuing the research in \cite{INo}, \cite{I3}, \cite{INUW} we show in this work that the stability is increasing when $k$ is growing in the presence of constant attenuation. For this problem, we could follow the arguments in \cite{NUW} to derive $k$-dependent stability estimates using only complex-valued geometrical optics solutions. However, in doing so, the constant associated with the H\"older part will grow exponentially in $k$, similar to the result in \cite{NUW}. To obtain polynomially growing constants as in \cite{I3} and \cite{INUW}, we use both complex- and real-valued geometrical optics solutions in our proof. In addition to geometrical optics solutions, we also need explicit sharp bounds on fundamental solutions of elliptic operators with parameter $k$. Let $\varepsilon$ be an operator norm of the difference of two Dirichlet-to-Neumann maps corresponding to different potentials. We give conditional estimates for difference of potentials by a function of $\varepsilon$ which goes to zero as $\varepsilon$ goes to zero. This function is the sum of the terms containing powers of $\varepsilon$ and of $-\log \varepsilon$, moreover the terms containing $ \log \varepsilon$ tend to zero (as powers of $k$) when $k\rightarrow \infty$.

In section 2 we state main results. In section 3 by using sharp bounds on regular fundamental solutions of some elliptic linear partial differential operators with complex coefficients containing large parameter $k$ we 
construct almost complex exponential solutions to the Schr\"odinger
equation with attenuation and give bounds on these solutions. In section 4 we similarly construct bounded almost complex exponential solutions
in the "low frequency zone" which is growing with increasing $k$. In section
5 we use these almost exponential solutions and the Fourier transform to prove our main results. In concluding section 6 we outline challenges and possible future research.

\section{ Main results}
Let $\Omega$ be a (bounded) domain in $\mathbb{R}^3$ with Lipschitz boundary. We  consider the Schr\"odinger equation
\begin{equation}
-\Delta u-k^2 u+ik bu+cu=0\;\mbox{in}\;\Omega
\label{schrodinger}
\end{equation}
with the Dirichlet boundary data
\begin{equation}
u=g\in H^{1/2}(\partial\Omega)\;\mbox{on}\;\partial\Omega.
\label{dirichlet}
\end{equation}
Assume that the attenuation coefficient $b$ is a non negative constant and the potential $c \in L_{\infty}(\Omega)$. Suppose that there exists a unique solution to \eqref{schrodinger}, \eqref{dirichlet}. Thus we can define the Dirichlet-to-Neumann map 
\begin{equation}
\Lambda_c g =\partial_{\nu}u\;\mbox{on}\;\partial \Omega.
\label{dtn}
\end{equation}
 By using potential theory one can show that $\Lambda_c-\Lambda_0$ is a continuous linear operator from $L^2(\partial\Omega)$ into $L^2(\partial\Omega)$. We denote its
operator norm by $||\Lambda_c-\Lambda_0||$.

We assume that $vol \Omega \leq 1$ and that $c$ is zero near $\partial\Omega$. Throughout we denote $C_0$ generic constants whose values may change from line to line. These constants do not depend  on $c$, $k$, or $\Omega$. They are only determined by our proofs. Generic constant $C(\Omega, M)$ might in addition depend on $\Omega, M$. We will use the norms $\|\cdot\|_p(\Omega)$ in the Lebesgue spaces $L^p(\Omega)$ and $\|\cdot\|_{(s)}(\Omega)$ in the Sobolev spaces $H^s(\Omega)$. 
\begin{theorem}{\rm(Low frequency/energy)}\label{main}
 Let
\begin{equation}
\|c_j\|_{\infty}(\Omega)\leq M_0,\;\|\nabla c_j\|_{\infty}( \Omega)\leq M_1,\;j=1,2,\;M=\sqrt{M_0^2+M_1^2},
\label{apriori}
\end{equation}
and $\varepsilon=\|\Lambda_{c_2}-\Lambda_{c_1}\|, E=-log \varepsilon$. Let
\begin{equation}
 k< E,\;\;2 \leq E.
\label{kE}
\end{equation}

  Then there is constant $C(\Omega, M)$ such that 
\begin{equation}
\|c_2-c_1\|_2^2(\Omega)\leq C(\Omega,M)e^{4b}\varepsilon (E^2+k^2)+
\frac{4M^2}{1+(\frac{E^2+k^2}{C(\Omega,M)})^{\frac{2}{3}}}.
\label{kbound}
\end{equation}
\end{theorem}

In the bound \eqref{kbound} the logarithmic component goes to zero as $k$ grows. Thus this bound can be viewed as an evidence of increasing stability in recoverying $c$ for larger frequencies/energies $k$.  In any event the second term, namely, $C\varepsilon k^2$, contributes only to (the best possible) Lipschitz stability, while logarithmic terms are decaying for larger $k$. 
While constants $C_0, C(\Omega, M)$ are hard to evaluate for general $\Omega$, it is very likely that when $\Omega$ is a ball or a cube  one can obtain relatively  simple explicit bounds on these constants. 

The factor $k$ in the second term of the bound \eqref{kbound} most likely is necessary. Indeed, one needs bounds on time derivatives  in the closely related inverse problems for the wave equation. For high frequencies/energies $k$, with additional assumptions on $E$ and $k$, one can derive a similar estimate.
\begin{theorem}{\rm(High frequency/energy)}\label{high}
\ Assume that \eqref{apriori} hold. Let 
 $\alpha>0, \beta \in (0,\frac{2}{3}]$. Let $k$ satisfy
\begin{equation}\label{Ek}
 E^{\alpha}\leq k^{\beta},\quad 2C_0^2 M^2 <k^2+2.
\end{equation}

Then there is constant $C(\Omega, M)$ such that 
\begin{equation}
\|c_2-c_1\|_{(0)}^2(\Omega)\leq 
C(\Omega,M)e^{4b} k^{3\beta}\varepsilon^2+
\frac{C(\Omega,M)}{(E^{\alpha}+k^{\beta})^2+C(\Omega,M)}.
\label{kEbound}
\end{equation}
\end{theorem}

We observe that for real-valued $c$ and $b=0$ the Dirichlet problem \eqref{schrodinger}, \eqref{dirichlet} might have eigenvalues $k$ when its solution fails to exist and be unique, so that the Dirichlet-to Neumann map is not well defined. Then one can consider instead the Neumann-to-Dirichlet map, or replace these maps by the Cauchy set with naturally defined norm.
   
\section{Almost complex exponential solutions}

We start with
\begin{lemma}\label{aces}
Let $\xi\in \mathbb{R}^3$ and 
\begin{equation}
\label{tau}
k^2 < \tau^2+\frac{|\xi|^2}{4},\; 8 C_0^2M^2<|\xi|^2+4\tau^2-2k^2+4
\end{equation}
for some $C_0>0$ determined later. Then there are solutions
\begin{equation}
\label{complexexp}
u(x;j)=e^{i x\cdot \zeta(j)}(1+v(x;j))
\end{equation}
to the equations
\begin{equation}
\label{ujsolution}
-\Delta u_j-k^2 u_j+ikb u_j+c_ju_j=0\;\mbox{in}\;\Omega\;
\end{equation} 
with
\begin{equation}
\label{xi}
\zeta(1)+\zeta(2)=\xi,\;
|\Im \zeta(j)|\le\sqrt{\frac{|\xi|^2}{4}+\tau^2-k^2+\frac{b}{2}k},
\end{equation}
\begin{equation}
\label{Ltails}
\|v( ;j)\|_{(1)}(\Omega)\leq \frac{2 C_0 M }{\sqrt{|\xi|^2+4\tau^2-2k^2+4}},
\end{equation}
and
\begin{equation}
\label{Htails}
\|v( ;j)\|_{(2)}(\Omega)\leq 2 C_0 M. 
\end{equation}
\end{lemma}
\noindent{\bf Proof}. Let $\xi\in{\mathbb R}^3, \xi\not=0$ and
$$
\xi=|\xi| e_1.
$$
We introduce
\begin{equation}\label{zeta}
\begin{aligned}
\zeta(1)&=\frac{|\xi|}{2}e_1+i(\frac{|\xi|^2}{4}+\tau^2-k^2+ikb)^{\frac{1}{2}}e_2+\tau e_3,\\
\zeta(2)&=\frac{|\xi|}{2}e_1-i(\frac{|\xi|^2}{4}+\tau^2-k^2+ikb)^{\frac{1}{2}}e_2-\tau e_3,
\end{aligned}
\end{equation}
where $e_1,e_2,e_3$ is an orthonormal basis in $\mathbb{R}^3$ and $\sqrt{z}$ is the principal branch of the square root function. Then (\ref{ujsolution}) holds if and only if
\begin{equation}
-\Delta v( ;j)-2i\zeta(j)\cdot\nabla v(;j)=c_j(1+v(;j))\;\mbox{in}\;\Omega.
\label{vjsolution}
\end{equation}
 
 To demonstrate the second bound \eqref{xi} we write  
 $(\frac{|\xi|^2}{4}+\tau^2-k^2+ikb)^{\frac{1}{2}} = X+iY, 0<X$. Squaring the
 both sides yields $X^2-Y^2= \frac{|\xi|^2}{4}+\tau^2-k^2, 2XY= kb$.
 Substituting $Y=\frac{kb}{2X}$ into the first equation and solving the resulting
 (quadratic) equation we obtain
 \[
 \begin{aligned}
 X &= \frac{\sqrt{|\xi|^2+4\tau^2-4k^2+\sqrt{(|\xi|^2+4\tau^2-4k^2)^2+16k^2b^2}}}
 {2\sqrt{2}}\\
 &\leq\frac{\sqrt{|\xi|^2+4\tau^2-4k^2+2k b}}{2},
 \end{aligned}
 \]
where we used the inequality $\sqrt{A^2+B^2}\leq A+B$ for any positive $A,B$. So we arrive at \eqref{xi}.

Let $P(\zeta;j)=\zeta\cdot\zeta+2\zeta(j)\cdot\zeta$. By known results \cite{Ho}, \cite{I1}, there is a regular fundamental solution
$E(j)$ of $P( ;j)$ such that for any linear partial differential operator $Q$ with constant coefficients
\begin{equation}
\|Q E(j)f\|_{(0)}(\Omega)\leq C_0 sup\frac{\tilde{Q}(\xi^*)}{\tilde{P}(\xi^*)}\|f\|_{(0)}(\Omega)\;(sup\; \mbox{over}\;\xi^*\in \mathbb{R}^3) 
\label{E}
\end{equation}
for any $f\in L^2(\Omega)$, where 
$$
\tilde{P}(\xi)=(\sum_{|\alpha|\leq 2}|\partial^{\alpha}_{\xi} P(\xi)|^2)^{\frac{1}{2}}.
$$
In our particular case, by letting $\zeta(j)=\xi(j)+i\eta(j), \xi(j), \eta(j) \in \mathbb{R}^3$, for any $\xi^*\in \mathbb{R}^3$, we have
$$
\tilde{P}^2(\xi^* ;j)=
(|\xi^*|^2+2\xi(j)\cdot\xi^*)^2+4(\eta(j)\cdot\xi^*)^2+
4(|\xi^*+\xi(j)|^2+|\eta(j)|^2)+12=
$$
$$
(|\xi^*|^2+2\xi^*\cdot \xi(j)+2)^2+4(\eta(j)\cdot\xi^*)^2+
4(|\xi(j)|^2+|\eta(j)|^2)+8\geq
$$
$$
4(|\xi(j)|^2+|\eta(j)|^2)+8=
$$
\begin{equation}\label{Lbound}
4\left(\frac{|\xi|^2}{4}+\tau^2+\sqrt{(\frac{|\xi|^2}{4}+\tau^2-k^2)^2+k^2b^2}\right)+8
\geq 2(|\xi|^2+4\tau^2-2k^2)+8,
\end{equation}
due to the choice of $\zeta(j)$ in (\ref{zeta}). On the other hand, observing $|\xi^*+\xi(j)|^2-|\xi(j)|^2=|\xi^*|^2+2\xi(j)\cdot\xi^*$ and using the elementary inequality $(A-B)^2+4A+2\geq 2(A+B)+1$, with  $A=|\xi^*+\xi(j)|^2, B=|\xi(j)|^2$, we obtain 
\begin{align}
\tilde{P}^2(\xi^* ;j)&\geq (|\xi^*+\xi(j)|^2-|\xi(j)|^2)^2+4|\xi^*+\xi(j)|^2+12\notag\\
&\ge 2(|\xi^*+\xi(j)|^2+|\xi(j)|^2)+11\geq |\xi^*|^2+1.\label{Hbound}
\end{align}

The regular fundamental solution in \cite{Ho} is a convolution operator, so it commutes with differentiations, and hence
from \eqref{E} it follows that
 \begin{equation}
\|Q E(j)f\|_{(1)}(\Omega)\leq C_0 sup\frac{\tilde{Q}(\xi^*)}{\tilde{P}(\xi^*)}
\|f\|_{(1)}(\Omega)\;(sup\; \mbox{over}\;\xi^*\in \mathbb{R}^3). 
\label{E1}
\end{equation}
Since $E(j)$ is a fundamental solution, any solution $v( ;j)$ to the equation
\begin{equation}
\label{fixed}
v( ;j)=E(j)(c_j(1+v( ;j)))\;\mbox{on}\; \Omega
\end{equation}
solves (\ref{vjsolution}). From (\ref{E1}) with $Q=1$ and (\ref{Lbound}) it follows that
$$
||E(j)f||_{(1)}(\Omega)\leq 
C_0(|\xi|^2+4\tau^2-2k^2+4)^{-\frac{1}{2}}||f||_{(1)}(\Omega),
$$
or
\begin{equation}
\label{theta}
||E(j)f||_{(1)}(\Omega)\leq \theta ||f||_{(1)}(\Omega)\;\mbox{with}\;
\theta = \frac{C_0}{\sqrt{|\xi|^2+4\tau^2-2k^2+4}}.
\end{equation}

Observe that
$$
\|c_j(1+v)\|_{(1)}(\Omega)\leq \|c_j\|_{(1)}(\Omega) + \|c_j v \|_{(1)}(\Omega)\leq M +\sqrt{2} M\|v\|_{(1)}(\Omega),
$$
where we used the apriori bound (\ref{apriori}). So the operator $F(v(;j))$ in the right side of (\ref{fixed}) maps the ball $B(\rho)=\{v:||v||_{(1)}(\Omega)\leq\rho\}$ into the ball
$B(\theta M  +\theta \sqrt{2} M \rho)$, and hence into $B(\rho)$ when
\begin{equation}
\label{into}
\theta M \leq (1-\theta \sqrt{2} M)\rho.
\end{equation}
 The second condition in (\ref{tau}) and  (\ref{theta}) imply that $\sqrt{2} \theta M \le \frac{1}{2}$, and hence (\ref{into}) holds with
\begin{equation}
\label{rho}
\rho=\frac{2 C_0 M }{\sqrt{|\xi|^2+4\tau^2-2k^2+4}}.
\end{equation}

Likewise, using (\ref{Hbound}) and (\ref{E}) with $Q(\xi^*)=\xi_k, k=1,2,3,$ from (\ref{fixed}) we have that for $v(;j)\in B(\rho)$
$$
||v( ;j)||_{(2)}(\Omega)\leq C_0||c_j (1+  v( :j))||_{(1)}(\Omega)\leq
C_0 M (1+\sqrt{2}||v( ;j)||_{(1)}(\Omega))\leq
$$
\begin{equation}
\label{vjH}
C_0 M (1+\frac{2\sqrt{2}C_0 M}{\sqrt{|\xi|^2+4\tau^2-2k^2+4}})\leq  2 C_0  M, 
\end{equation}
due to the  bound (\ref{Ltails}) given by (\ref{rho}) and the second
inequality \eqref{tau}.  Therefore the operator $F$ is continuous from $H^1(\Omega)$ into $H^2(\Omega)$ and therefore compact from $H^1(\Omega)$ into itself.
Now, this operator maps convex closed set $B(\rho)\subset H^1(\Omega)$ into itself 
and is compact, hence by Schauder-Tikhonov Theorem it has a 
fixed point $v( ;j)\in B(\rho)$. In view of (\ref{rho}) we have the bound (\ref{Ltails}) and due to (\ref{vjH}) we have the bound (\ref{Htails}). The proof is complete.\eproof

In the following lemma we will derive boundary estimates of almost complex exponential solutions constructed in Lemma~\ref{tau}.
\begin{lemma}
Let $u(;j)$ be the solutions (\ref{complexexp}) to the Schr\"odinger equations 
$-\Delta u_j-k^2 u_j+ibk u_j+c_ju_j=0\;\mbox{in}\;\Omega$ from Lemma 3.1. Then one has
\begin{equation}
\label{uboundary}
\|u( ;j)\|_{(0)}(\partial\Omega)\leq(|\partial\Omega|^{1/2}+ C(\Omega) )
e^{\sqrt{\frac{|\xi |^2}{4}+\tau^2-k^2+\frac{1}{2}bk}}.
\end{equation}
\end{lemma}

\noindent{\bf Proof}.
Let $e(x;j)=e^{i x\cdot\zeta(j)}$. Obviously,
$$
||e( ;j)||_{(0)}(\partial \Omega)\leq |\partial\Omega|^{\frac{1}{2}}
e^{|\Im \zeta(j)|}.
$$
Moreover, from trace theorems for Sobolev spaces
$$
\|v( ;j)\|_{(0)}(\partial \Omega)
\leq C(\Omega) \|v( ;j)\|_{(1)}(\Omega)
$$
$$
\leq C(\Omega) \frac{C_0 M }{\sqrt{|\xi|^2+4\tau^2-2k^2+4}}\leq\frac{1}{2}C(\Omega)
$$
by \eqref{Ltails} and \eqref{tau}. Hence, we have  
$$
\|u( ;j)||_{(0)}(\partial \Omega)\leq \|e( ;j)\|_{(0)}(\partial\Omega)+
e^{|\Im \zeta(j)|}\|v( ;j)\|_{(0)}(\partial\Omega)\leq (|\partial\Omega|^{1/2}+C(\Omega))e^{|\Im \zeta( ;j)|}
$$
and \eqref{uboundary} follows immediately from the second estimate of \eqref{xi}. \eproof

\section{Bounded complex exponential solutions}

In this section we construct bounded complex exponential solutions. As in the previous section, we begin with
\begin{lemma}\label{ares}
Let $\xi\in \mathbb{R}^3$ and 
\begin{equation}
\label{tauR}
|\xi|^2 \leq 3 k^2,\; 2C_0^2M^2< k^2+2.
\end{equation}
Then there are solutions
\begin{equation}
\label{realexp}
u(x;j)=e^{i x\cdot \zeta(j)}(1+v(x;j))
\end{equation}
to the equations
\begin{equation}
\label{ujsolutionR}
-\Delta u_j-k^2 u_j+ikb u_j+c_ju_j=0\;\mbox{in}\;\Omega,\;
\end{equation} 
with
\begin{equation}
\label{xiR}
\zeta(1)+\zeta(2)=\xi,\;
|\Im \zeta(j)| \leq b ,
\end{equation}
\begin{equation}
\label{LtailsR}
\|v( ;j)\|_{(1)}(\Omega)\leq \frac{C_0 M }{\sqrt{k^2+2}}
\end{equation}
and
\begin{equation}
\label{HtailsR}
\|v( ;j)\|_{(2)}(\Omega)\leq 2C_0 M. 
\end{equation}
\end{lemma}

\noindent{\bf Proof}. Let $\xi\in{\mathbb R}^3, \xi\not=0$ and
$$
\xi=|\xi| e_1.
$$
We introduce
$$
\zeta(1)=\frac{|\xi|}{2}e_1+(k^2-\frac{|\xi|^2}{4}-ikb)^{\frac{1}{2}}e_2,
$$
\begin{equation}
\label{zetaR}
\zeta(2)=\frac{|\xi|}{2}e_1-
(k^2-\frac{|\xi|^2}{4}-ikb)^{\frac{1}{2}}e_2.
\end{equation}
Then (\ref{ujsolutionR}) holds if and only if
\begin{equation}
-\Delta v( ;j)-2i\zeta(j)\cdot\nabla v(;j)=c_j(1+v(;j))\;\mbox{in}\;\Omega.
\label{vjsolutionR}
\end{equation}

To demonstrate the second bound of \eqref{xiR}, we write $(k^2-\frac{|\xi|^2}{4}-ikb)^{1/2}=X+iY$ with $0<X$. Squaring the both sides implies
$X^2-Y^2=k^2-\frac{|\xi|^2}{4}$, $2XY=-kb$. Substituting $X=-\frac{kb}{2Y}$ into the first equation and solving the resulting (quadratic) equation, we obtain
\[
\begin{aligned}
|Y|=\frac{kb}{2X}&=\frac{kb}{\sqrt{2}\sqrt{(k^2-\frac{|\xi|^2}{4})+\sqrt{(k^2-\frac{|\xi|^2}{4})^2+k^2b^2}}}\\
&\leq\frac{kb}{2\sqrt{k^2-\frac{|\xi|^2}{4}}}\le b
\end{aligned}
\]
in view of the first condition of \eqref{tauR}. 

As before, let $P(\zeta;j)=\zeta\cdot\zeta+2\zeta(j)\cdot\zeta$
and  $\zeta(j)=\xi(j)+i\eta(j), \xi(j), \eta(j) \in {\mathbb R}^3$. Observe that
$$
|\xi(j)|^2 + |\eta(j)|^2 = \frac{|\xi|^2}{4}+
\sqrt{(k^2 - \frac{|\xi|^2}{4})^2+k^2b^2}\geq   
\frac{|\xi|^2}{4}+\sqrt{(k^2 - \frac{|\xi|^2}{4})^2} = k^2   
$$
and using \eqref{Lbound} we yield
\begin{equation}
\label{LboundR}
\tilde{P}^2(\xi^* ;j)\geq 4(|\xi(j)|^2+|\eta(j)|^2)+8 
\geq 4 k^2+8.
\end{equation}

Similarly, any solution $v( ;j)$ to the equation
\begin{equation}
\label{fixedR}
v( ;j)=E(j)(c_j(1+v( ;j)))\;\mbox{on}\; \Omega
\end{equation}
satisfies (\ref{vjsolutionR}). From (\ref{E1}) with $Q=1$ and (\ref{LboundR}) it follows that
$$
||E(j)f||_{(1)}(\Omega)\leq C_0(4k^2+8)^{-\frac{1}{2}}||f||_{(1)}(\Omega),
$$
or
\begin{equation}
\label{thetaR}
||E(j)f||_{(1)}(\Omega)\leq \theta ||f||_{(1)}(\Omega),\;\mbox{where}\;
\theta = \frac{C_0}{\sqrt{4k^2+8}}.
\end{equation}

As in the proof of Lemma~\ref{aces} the operator $F(v(;j))$ in the right side of (\ref{fixedR}) maps the ball $B(\rho)=\{v:||v||_{(1)}(\Omega)\leq\rho\}$   into itself when $ \theta M \leq (1-\theta \sqrt{2} M)\rho $, which holds with
\begin{equation}
\label{rhoR}
\rho=\frac{C_0 M }{\sqrt{k^2+2}}.
\end{equation}
because $1 - \theta\sqrt{2} M> \frac{1}{2}$ due to the second condition
\eqref{tauR}. 

Repeating similar arguments in Lemma~\ref{aces}, the operator $F$ is continuous from $H^1(\Omega)$ into $H^2(\Omega)$ and therefore compact from $H^1(\Omega)$ into itself.
Now this operator maps convex closed set $B(\rho)\subset H^1(\Omega)$ into itself 
and is compact, hence by Schauder-Tikhonov Theorem it has a fixed point $v( ;j)\in B(\rho)$. Due to (\ref{rhoR}) we have the bound (\ref{LtailsR}) and due to (\ref{vjH}) we have the bound (\ref{HtailsR}).\eproof

We also estimate the boundary values of the almost real exponential solution $u(;j)$.
\begin{lemma}
Let $u(;j)$ be the solutions (\ref{realexp}) to the Schr\"odinger equations 
$-\Delta u_j-k^2 u_j+c_ju_j=0\;\mbox{in}\;\Omega$ from Lemma~\ref{ares}. Then one has
\begin{equation}
\label{uboundaryR}
\|u( ;j)\|_{(0)}(\partial\Omega)\leq (|\partial\Omega|^{\frac{1}{2}}+C(\Omega)) e^{b}.
\end{equation}
\end{lemma}
\noindent{\bf Proof}. Let $e(x;j)=e^{i x\cdot\zeta(j)}$. Obviously,
$$
||e( ;j)||_{(0)}(\partial \Omega)\leq |\partial\Omega |^{\frac{1}{2}}
e^{|\Im \zeta(j)|}.
$$

Moreover, from trace theorems for Sobolev spaces
$$
\|v( ;j)\|_{(0)}(\partial\Omega)
\leq C(\Omega) \|v( ;j)\|_{(1)}(\Omega)
$$
$$
\leq C(\Omega) \frac{C_0 M }{\sqrt{k^2+2}}\leq C(\Omega)
$$
by \eqref{LtailsR} and \eqref{tauR}. Hence we have that
\begin{align*}
\|u( ;j)\|_{(0)}(\partial \Omega)&\leq \|e( ;j)\|_{(0)}(\partial\Omega)+
e^{|\Im \zeta(j)|}\|v( ;j)\|_{(0)}(\partial\Omega)\\
&\leq
|\partial\Omega|^{\frac{1}{2}} e^{|\Im \zeta(j)|}+ 
C(\Omega)e^{|\Im \zeta(j)|}
\end{align*}
and \eqref{uboundary} holds from the second bound of \eqref{xiR}.\eproof

\section{Proofs of stability estimates}

The following standard orthogonality result (see, for example, \cite{A}, \cite{I1}) follows by simple application of the Green's formula.

\begin{lemma}
\begin{equation}
\label{Green}
\int_{\Omega}(c_1-c_2) u_1 u_2=\int_{\partial\Omega}((\Lambda_{c_2}-\Lambda_{c_1})u_1)u_2
\end{equation}
for all functions $u_1, u_2 \in H^{1}(\Omega)$, solving the Schr\"odinger equations
$$
-\Delta u_1-k^2u_1 +ibu_1+c_1u_1=0\;\mbox{in}\;\Omega
$$ 
and
$$
-\Delta u_2-k^2 u_2+ ib u_2+c_2u_2=0\;\mbox{in}\;\Omega.
$$ 
\end{lemma}

\noindent{\bf Proof of Theorem~\ref{main}}

Substituting almost complex exponential solutions (\ref{complexexp}) into the identity (\ref{Green}) and observing that
$$
u(x;1)u(x;2)= e^{ix\cdot\xi}(1+v(x;1)+v(x;2)+v(x;1)v(x;2)),
$$
we obtain for the Fourier transform $(\hat{c}_2-\hat{c}_1)(\xi)$ of $c_2-c_1$
$$
|(\hat{c}_2-\hat{c}_1)(\xi)|=|\int_{\Omega}(c_1-c_2)(x)e^{i\xi\cdot x} dx|\leq
$$
\begin{equation}
\label{bound1}
 \int_{\Omega}|c_2-c_1|(|v(;1)|+|v(;2)|+|v(;1)||v(;2)|)+
\varepsilon \|u( ;1)\|_{(0)}(\partial \Omega)\|u( ;2)\|_{(0)}(\partial \Omega).
\end{equation}

With the help of \eqref{Htails}, Sobolev Embedding Theorems imply
\begin{equation}\label{embedding}
\|v( ;2)\|_{\infty}(\Omega)\leq C_e(\Omega)\|v( ;2)\|_{(2)}(\Omega)\leq
2C_e(\Omega)C_0 M.
\end{equation}
Therefore,
\begin{align}\label{bound0}
&\int_{\Omega}|c_2-c_1|(|v(;1)|+|v(;2)|+|v(;1)||v(;2)|)\notag\\
&=\int_{\Omega}|c_2-c_1||v(;1)|+\int_{\Omega}|c_2-c_1||v(;2)|+\int_{\Omega}|c_2-c_1||v(;1)||v(;2)|\notag\\
&\leq\|c_2-c_1\|_{(0)}(\Omega)\|v(;1)\|_{(0)}(\Omega)+\|c_2-c_1\|_{(0)}(\Omega)\|v(;2)\|_{(0)}(\Omega)\notag\\
&\quad+\|c_2-c_1\|_{(0)}(\Omega)2C_e(\Omega)C_0 M \|v(;1)\|_{(0)}(\Omega)\notag\\
&\leq\|c_2-c_1\|_{(0)}(\Omega)4(1+C_e(\Omega)C_0 M)\frac{C_0 M }{\sqrt{|\xi|^2+4\tau^2-2k^2+4}},
\end{align}
where we used the Cauchy-Schwarz inequality, \eqref{embedding}, and, in the last inequality, \eqref{Ltails}. From  \eqref{bound1}, \eqref{bound0}, and \eqref{uboundary},  we have that
$$
|(\hat{c}_2-\hat{c}_1)(\xi)|\leq \|c_2-c_1\|_{(0)}(\Omega)4(1+C_e(\Omega)C_0 M)
\frac{C_0 M }{\sqrt{|\xi|^2+4\tau^2-2k^2+4}}+
$$
$$
\varepsilon C^2(\Omega) e^{\sqrt{|\xi|^2+4\tau^2-4k^2+2bk}},
$$
provided \eqref{tau} is satisfied.

So  by using the Parseval identity and polar coordinates we conclude that
\begin{align}\label{bound2}
\|c_2-c_1\|_{(0)}^2(\Omega)&=\int|\hat{c}_2-\hat{c}_1|^2(\xi) d\xi\leq
 \int_{|\xi|<\rho}|\hat{c}_2-\hat{c}_1|^2(\xi) d\xi+\int_{\rho<|\xi|}|\hat{c}_2-\hat{c}_1|^2(\xi) d\xi\notag\\
 &\le\|c_2-c_1\|^2_{(0)}(\Omega)C_3 \int_0^{\rho}\frac{r^2}{r^2+4\tau^2-2k^2+4} dr\notag\\
 &\quad+
8\pi\varepsilon^2C^4(\Omega)\int_0^{\rho}e^{2\sqrt{r^2+4\tau^2-4k^2+2bk}} r^2 dr+
\int_{\rho<|\xi|}|\hat{c}_2-\hat{c}_1|^2(\xi) d\xi.
\end{align}
where $C_3(\Omega,M)=32\pi (1+C_e(\Omega)C_0 M)^2C_0^2M^2$. Let
\begin{equation}
\label{choice}
4\tau^2=4k^2-r^2+(\frac{E}{2})^2-2kb+4b^2+8C^2_0M^2,
\end{equation}
then the second bound of \eqref{tau} is clearly satisfied, and the first bound of \eqref{tau} holds provided $k < E$, as guaranteed by
\eqref{kE}. Now we choose
\begin{equation}
\label{choicerho}
\rho= R^{-1}(E^2+k^2)^{\frac{1}{3}}.
\end{equation}
with some $R>2$ to be selected later on. Then the right side of 
\eqref{choice} is positive and hence our choice of $\tau$ is possible.
Indeed, due to \eqref{choicerho} it suffices to show that
$$
\frac{(E^2+k^2)^{\frac{2}{3}}}{4}+2bk< 4k^2+\frac{E^2}{4}+4b^2+8C^2_0M^2. 
$$
We recall the elementary inequality $(A+B)^{\lambda}\leq A^{\lambda}+ 
B^{\lambda}$ for any positive numbers $A,B$ and $\lambda\in (0,1)$. Using this
inequality with $\lambda =\frac{2}{3}$ and $A = E^2, B=k^2$ we conclude
that the needed inequality follows from
$$
E^{\frac{4}{3}}+k^{\frac{4}{3}}+8bk< 16k^2+ E^2+ 16 b^2+32 C^2_0M^2. 
$$
This inequality in turn follows from three inequalities
$$
E^{\frac{4}{3}}< a E^2,\ k^{4/3}<11k^2+(1-a)E^2+32C_0^2M^2,\ \text{and}\ 8bk<5k^2+16 b^2,
$$
where $a=2^{-2/3}$. The first and second ones hold since we assumed $2<E$ and the third one is obvious.

The second condition \eqref{tau} follows from our choice of 
$\tau$ in \eqref{choice}. Indeed, due to \eqref{choice} this condition
becomes
$$
8 C_0^2M^2< 2k^2-2bk+4b^2 +\frac{E^2}{4}+ 8 C_0^2 M^2+4
$$
which is obvious, since $0 \leq 2k^2-2bk+ 4b^2 $.

Due to \eqref{choice}
$$
\int_0^{\rho}\frac{r^2}{r^2+4\tau^2-2k^2+4} dr\leq
\int_0^{\rho}\frac{r^2}{2k^2-2bk+b^2+(\frac{E}{2})^2+4} dr=
$$
\begin{equation}
\label{integral1}
\frac{1}{2k^2-2bk+b^2+(\frac E2)^2+4}\frac{\rho^3}{3}=\frac{E^2+k^2}{3(k^2+(k-b)^2+(\frac{E}{2})^2+4)R^3}\leq \frac{2}{R^3},
\end{equation}
and
$$
\varepsilon^2 \int_0^{\rho}e^{2\sqrt{r^2+4\tau^2-4k^2+2bk}} r^2 dr\leq
\varepsilon^2 \int_0^{\rho}e^{\sqrt{E^2+16b^2+ 32C_0^2M^2}} r^2 dr\leq
$$ 
\begin{equation}
\label{integral2}
\varepsilon^2 \int_0^{\rho}e^{E+6C_0 M+4b} r^2 dr=\varepsilon e^{6C_0 M+4b}\frac{\rho^3}{3}.
\end{equation}
From \eqref{apriori}, we have
$$
\int_{\rho<|\xi|}|\hat{c}_2-\hat{c}_1|^2(\xi) d\xi\leq \frac{4M^2}{\rho^2+1}.
$$
Thus we obtain from (\ref{bound2}) that
$$
\|c_2-c_1\|_{(0)}^2(\Omega)\leq 
 \|c_2-c_1\|_{(0)}^2(\Omega)\frac{2C_3}{R^3}+
$$
$$
8\pi C(\Omega)^4\frac{1}{3R^3}e^{6C_0M+4b}\varepsilon(E^2+k^2)+
\frac{4 M^2}{1+R^{-2}(E^2+k^2)^{\frac{2}{3}}}.
$$
Choosing $R=\max\{ (4C_3)^{\frac{1}{3}},2\}$ we will absorb the first term on the right side by the left side.  The proof is complete.\eproof

\noindent{\bf Proof of Theorem~\ref{high}}.

Using the bounded almost exponential solutions (\ref{realexp}) in the identity (\ref{Green}) and \eqref{embedding} as above  we obtain that
$$
\begin{aligned}
|(\hat{c}_2-\hat{c}_1)(\xi)|&\leq 
\|c_2-c_1\|_{(0)}(\Omega)\left(2C_e(\Omega)M+1)\|v(;1)\|_{(0)}(\Omega)+
\|v(;2)\|_{(0)}(\Omega)\right)+\\
&\qquad\varepsilon \|u( ;1)\|_{(0)}(\partial \Omega)\|u( ;2)\|_{(0)}(\partial \Omega)\\
&\leq\|c_2-c_1\|_{(0)}(\Omega)(2+2C_e(\Omega)C_0M)
\frac{C_0 M }{\sqrt{k^2+2}}+\varepsilon 2C(\Omega)^2e^{2b} ,
\end{aligned}
$$
where we have utilized \eqref{LtailsR} and \eqref{uboundaryR}.

As above, using the Parseval identity and polar coordinates we conclude that
$$
\|c_2-c_1\|_{(0)}^2(\Omega)\leq
 \int_{|\xi|<\rho}|\hat{c}_2-\hat{c}_1|^2(\xi) d\xi+\int_{\rho<|\xi|}|\hat{c}_2-\hat{c}_1|^2(\xi) d\xi\leq
$$
\begin{equation}
\label{bound2R}
\|c_2-c_1\|^2_{(0)}(\Omega)C_3 \int_0^{\rho}\frac{r^2}{k^2+2} dr+
8\pi\varepsilon^2C^4(\Omega)e^{4b} \int_0^{\rho} r^2 dr+
\int_{\rho<|\xi|}|\hat{c}_2-\hat{c}_1|^2(\xi) d\xi.
\end{equation}
where $C_3=32\pi (1+C_e(\Omega)C_0M)^2C_0^2M^2$. Hence we obtain from (\ref{bound2R})
$$
\|c_2-c_1\|_{(0)}^2(\Omega)\leq \|c_2-c_1\|_{(0)}^2(\Omega)\frac{C_3\rho^3}{3(k^2+2)}+
\frac{8\pi}{3}C(\Omega)^4e^{4b} \rho^3 \varepsilon^2 +\frac{4 M^2}{\rho^2+1}.
$$
Choosing $\rho = (2C_3^{\frac 13})^{-1}(E^{\alpha}+k^{\beta}) $, we can absorb the first term on the right side by the left side and obtain the bound \eqref{kEbound}. It only remains to show that with this choice the first bound of \eqref{tauR} is satisfied.  To this end, in view of \eqref{Ek}, we deduce that
\[
\rho^2=(2C_3^{\frac 13})^{-2}(E^{\alpha}+k^{\beta})^2\le C_3^{-\frac 23}k^{2\beta}\le(\frac{1}{32\pi C_0^2M^2k})^{\frac 23}k^2\le 3k^2.
\]
The proof is complete.\eproof

\section{Conclusion}

In this paper we showed that stability of recovery of the Schr\"odinger potential is increasing in presence of constant attenuation which is a feature of most applied problems.   There is a belief that stability in the continuation and inverse problems always grows with frequency $k$. As shown in \cite{J}, in general stability of the continuation for the Helmholtz equation might deteriorate.  In  \cite{HI}, \cite{I2} it was shown that the stability of the continuation is improving under some (convexity  type) conditions. 
In \cite{IK} it was demonstrated that in some cases these convexity conditions can be relaxed.

Now we outline  possible future developments and challenges.

 We believe that constants $C_0$ in the stability estimates may be evaluated more explicitly by using periodic Faddeev type solutions \cite{H} or more suitable  regular fundamental solutions based on the work of H\"ormander
 in 1955 \cite{Ho}. We expect to obtain explicit bounds  at least when $\Omega$ is the unit ball.  One expects to obtain improving
 stability results when $c$ is not necessarily zero near $\partial\Omega$
 by using methods of \cite{SU1} or of singular solutions \cite{I1} for boundary reconstruction. It would be interesting to use the new ideas
 in \cite{Bu} to handle the two-dimensional case.
 It is not clear now how to get better stability
 for both $c$ and constant $b$ or for $b$ close to a constant.
 Most likely, one will need to use the Fourier Integral Operators instead of the Fourier Transform. 
 
 There is a need in an additional  numerical evidence of increasing stability in the important inverse problem we considered in this paper, as well as in similar inverse medium problems  studied in particular in \cite{B}, \cite{Na}. In our view, even numerical results for the linearized problem (like Born approximation) would be convincing and  interesting.
 
 We hope to demonstrate the increasing stability for hard or transparent (convex) obstacles from the Dirichlet-to-Neumann map by combining the results of \cite{I2}, \cite{IK} with the methods of \cite{AD}. So far increasing
 stability of reconstruction of obstacles was observed numerically, but there are no analytic results explaining it and suggesting better numerical methods.

It is still an  open question whether logarithmic stability of recovery of near  field from far field pattern  \cite{Bus}, \cite{I1}, section 6.1 is improving with growing frequency. To show it one can try to adjust the methods of \cite{IK} to handle Hankel functions.
 
 Probably, it will be difficult to show increasing stability for the coefficient $a_0$ in the equation $ (-\Delta -k^2 a_0^2(x))u=0$.  At present, there are  only some preliminary results (in the plane case) 
 \cite{P} under some non trapping conditions on $a_0$ and bounds with constants exponentially growing with $k$ \cite{NUW}. We observe that the conductivity equation $div(a\nabla u)+k^2 u=0$ can be transformed into this Helmholtz type equation with $a_0=a^{-\frac{1}{2}}, c = a^{-\frac{1}{2}}\Delta a^{\frac{1}{2}}$, so increasing stability for the conductivity equation remains a challenge.

\vspace{6mm} \noindent{\Large{\bf Acknowledgments}}

\medskip
This research is supported in part by the Emylou Keith and Betty Dutcher Distinguished Professorship and the NSF grant DMS 10-08902. The second author is partially supported by the National Science Council of Taiwan.

\end{document}